%
%
%
\magnification=1200
\input amstex
\documentstyle{amsppt}
\topmatter
\title
Quadrature Formulas Based on \\ Rational Interpolation
\endtitle
\author
Walter Van Assche \\ Ingrid Vanherwegen
\endauthor
\affil
Katholieke Universiteit Leuven
\endaffil
\address
Department of Mathematics, Katholieke Universiteit Leuven,
Celestijnenlaan 200\,B, B-3001 Heverlee (BELGIUM)
\endaddress
\email
fgaee03\@cc1.kuleuven.ac.be
\endemail
\thanks
The first author is a Research Associate of the Belgian National Fund
for Scientific Research
\endthanks
\leftheadtext{Walter Van Assche -- Ingrid Vanherwegen}
\rightheadtext{Quadrature Formulas}

\abstract
We consider quadrature formulas based on interpolation using the basis
functions
$1/(1+t_kx)$ $(k=1,2,3,\ldots)$ on $[-1,1]$, where $t_k$ are parameters on the
interval $(-1,1)$. We investigate  two types of quadratures:
quadrature formulas of maximum accuracy
which correctly integrate as
many basis functions as possible (Gaussian quadrature), and
quadrature formulas whose
nodes are the zeros of the orthogonal functions  obtained by
orthogonalizing
the system of basis functions (orthogonal quadrature). We show that both
approaches involve orthogonal polynomials with modified (or varying) weights
which depend on the number of quadrature nodes.
The asymptotic distribution of the nodes is obtained as well as various
interlacing properties and monotonicity results for the nodes.
\endabstract

\endtopmatter

\document
\head 1. Introduction \endhead

Suppose that we want to compute the integral
$$   \int_{-1}^1 f(x) w(x) \, dx $$
with $w(x)$ a positive and integrable weight function on $[-1,1]$; then one
can use an interpolatory quadrature formula by approximating $f$
using interpolation. In this paper we assume that $w(x) > 0$ almost everywhere
in $[-1,1]$ and use interpolation by means of the basis functions
$$ \frac1{1+t_1x}, \frac1{1+t_2x}, \ldots, \frac1{1+t_nx}, \ldots,   \tag 1.1 $$
where $t_1,t_2,t_3,\ldots$ are parameters belonging to $(-1,1)$. We consider two
approaches for choosing the interpolation nodes (quadrature nodes). The first
consists of choosing the $n$ nodes to be the zeros of the orthogonal
function $r_{n+1}$ obtained by orthogonalizing the system \thetag{1.1} using the
inner product
$$   (f,g) = \int_{-1}^1 f(x)g(x) w(x) \, dx . $$
We refer to this as orthogonal quadrature. The second approach that we
investigate is to choose the nodes in such a way that the quadrature formula is
correct for as many basis functions as possible. This gives a
quadrature with maximum accuracy, and we refer to it as Gaussian quadrature.
It is well known that both approaches are the
same when one is using interpolation by polynomials, but these choices result in
different quadrature formulas when one uses interpolation with the basis
functions \thetag{1.1}.

\subhead 1.1 Orthogonal functions \endsubhead
It is well known that the zeros of the $n$th-degree orthogonal polynomial
with respect to a positive weight function $w(x)$ on $[-1,1]$ are all real,
simple and in the interval $(-1,1)$, and that they separate the zeros of the
orthogonal polynomial of degree $n+1$ \cite{26, p.\ 44--46}.
These properties can be generalized to other orthogonal systems of functions
$\{ \phi_i: i=1,2,\ldots\}$, where orthogonality is with respect to
the weight function $w(x)$ on $[-1,1]$, provided they satisfy the {\it Haar
condition}. The sequence of functions
$\{ \phi_i: i=1,\ldots,n\}$ satisfies the Haar condition on $[-1,1]$ if and
only if for every choice of $n$ points $x_i \in [-1,1]$ $(i=1,\ldots,n)$
with $x_i \neq x_j$ whenever $i \neq j$ one has $\det_{1 \leq i,j \leq n}
\left(\phi_i(x_j)\right) \neq 0$ (see, e.g., \cite{14}).
 This condition is equivalent to saying that
each linear combination
$\phi = \sum_{i=1}^n a_i\phi_i$ has at most $n-1$ zeros in $[-1,1]$.
A system $\{ \phi_i: i=1,\ldots,n\}$ is an {\it extended Haar system} on
$[-1,1]$ if and only if each linear combination
$\phi = \sum_{i=1}^n a_i\phi_i$ has at most
$n-1$ zeros in $[-1,1]$, counting multiplicities.
Pinkus and Ziegler \cite{25} and Videnski\u\i\  \cite{28} have investigated
several
properties of systems of {\it orthogonal} functions
$\{ \phi_i: i=1,2,\ldots,n+1\}$ satisfying the Haar condition on $[-1,1]$.
They have shown that $\phi_{n+1}$ has precisely $n$ zeros in $(-1,1)$, the
multiplicity of each zero is odd and the zeros of $\phi_n$ and $\phi_{n+1}$
are strictly interlacing.

Suppose that $t_i \in (-1,1)$ $(i=1,2,\ldots)$ and that $t_i \neq t_j$ whenever
$i \neq j$. Then the system
$\left\{ \frac1{1+t_1x}, \frac1{1+t_2x}, \ldots, \frac1{1+t_{n+1}x}\right\}$
is an extended Haar system on $[-1,1]$ because
$$      \phi(x) = \sum_{i=1}^{n+1} a_i\phi_i(x) =
\frac{p_n(x)}{\prod_{i=1}^{n+1} (1+t_ix)}\ , $$
where $p_n(x)$ is a polynomial of degree at most $n$, and hence
$\phi$ has at most $n$ zeros on $[-1,1]$, counting multiplicities.
If we apply the Gram-Schmidt procedure to
$\left\{ \frac1{1+t_1x}, \frac1{1+t_2x}, \ldots \right\}$
then we obtain a sequence of orthogonal functions on $[-1,1]$ with
respect to the weight function $w(x)$. We denote these orthogonal functions by
$$ r_{n+1}(t_1,\ldots,t_{n+1};x) = \frac{1}{1+t_{n+1}x} + \sum_{i=1}^{n}
\frac{c_i}{1+t_ix}, $$
and, as mentioned before, this rational function has precisely $n$ zeros
in $(-1,1)$.
It is well known that
$$    \min_{(a_1,\ldots,a_n) \in {\Bbb R}^n}
\int_{-1}^1 \left| \frac{1}{1+t_{n+1}x} - \sum_{i=1}^{n}
\frac{a_i}{1+t_ix} \right|^2 w(x)\, dx  $$
is obtained when
$$  \frac{1}{1+t_{n+1}x} - \sum_{i=1}^{n} \frac{a_i}{1+t_ix}
    = r_{n+1}(t_1,\ldots,t_n;x) ,  $$
so that $r_{n+1}(t_1,\ldots,t_n;x)$ is the error function for the best
$L_2(w)$-approximation of $\frac{1}{1+t_{n+1}x}$ by linear combinations of
$\left\{ \frac1{1+t_1x}, \frac1{1+t_2x}, \ldots, \frac{1}{1+t_nx} \right\}$.

\subhead 1.2 Gaussian quadrature \endsubhead
Suppose the system of functions $\{ \phi_i : i=1,2,\ldots\}$ is an
{\it extended Markov system} on $[-1,1]$, i.e. $\{ \phi_i : i=1,2,\ldots,k\}$ is
an extended Haar system on $[-1,1]$ for each $k \in \Bbb N$.
If one considers a quadrature formula with $n$ nodes which is of
maximum accuracy, then one wants to
integrate the functions $\phi_i$ $(1 \leq i \leq m)$
correctly for $m$ as large as possible. If $\phi_i(x) = x^{i-1}$, then
$m$ turns out to be $2n$, and the nodes are the zeros of the orthogonal
polynomial of degree $n$ for the weight function $w(x)$ (Gauss-Jacobi
quadrature). Krein \cite{15} showed that when $\{ \phi_i : i=1,2,\ldots,2n\}$
is a Haar system on $[-1,1]$, then there exists a unique set of $n$ nodes $x_i
\in (-1,1)$ and $n$ strictly positive weights $\beta_i$ such that
for $k=1,2,\ldots,2n$
$$   \int_{-1}^1 \phi_k(x)w(x)\,dx = \sum_{i=1}^n \beta_i \phi_k(x_i),
     \tag 1.2 $$
and this gives maximum accuracy. Karlin and Pinkus \cite{12} have generalized
this by allowing nodes with multiplicities and there are further generalizations
by Barrow \cite{2} and Bojanov, Braess and Dyn \cite{3}.

The system of rational functions \thetag{1.1} is an extended Markov system
on $[-1,1]$, hence the properties mentioned above are all valid. In particular
it follows that the nodes for Gaussian quadrature are all in $(-1,1)$
and the corresponding quadrature weights are all strictly positive.

\head 2. Connection with orthogonal polynomials \endhead

In this section we will make the simple observation that the nodes and weights
for both, orthogonal quadrature and Gaussian quadrature, based on the rational
functions \thetag{1.1} are closely
related to the zeros and Christoffel numbers for orthogonal polynomials
on $[-1,1]$ with a weight function that depends on the number of nodes.
This observation is not new (see e.g., L\'opez and Ill\'an \cite{18} \cite{19}
for Gaussian quadrature), but this interaction makes it possible to
use results from the theory of orthogonal polynomials to obtain useful
properties of the nodes and weights for  quadrature based on rational
interpolation.

\proclaim{Theorem 1}
Let $\{\phi_k\}$ be the set of basis functions \thetag{1.1}.
\roster
\item Suppose all the $t_i$ are different; then the $n$ nodes for orthogonal
quadrature are the zeros $x_{j,n}$ $(1 \leq j \leq n)$
of the orthogonal polynomial $p_{n}(x)$ with respect
to the weight function
$$   w_n(x)=  \frac{w(x)}{\pi_n(x)\pi_{n+1}(x)}, \qquad -1 < x < 1, \tag 2.1 $$
where
$$    \pi_n(x) = \prod_{i=1}^n (1+t_ix) . \tag 2.2 $$
If $t_1=0$ then \thetag{1.2} holds for $k=1,2,\ldots,n+1$ with $x_j=x_{j,n}$
and
$$    \beta_j(t_1,\ldots,t_{n+1}) = \pi_{n+1}(x_{j,n})\pi_n(x_{j,n})
\lambda_{j,n},   $$
where $\lambda_{j,n}$ are the Christoffel numbers corresponding to the
orthogonal polynomial $p_n(x)$
for the weight function \thetag{2.1}.
\item Suppose all the $t_i$ are different; then the nodes for Gaussian
quadrature are the zeros $x_{j,n}$ $(1 \leq j \leq n)$
of the orthogonal polynomial $p_{n}(x)$ with respect
to the weight function
$$   w_n(x)=  \frac{w(x)}{\pi_{2n}(x)}, \qquad -1 < x < 1, \tag 2.3 $$
and \thetag{1.2} holds for $k=1,2,\ldots,2n$ with $x_j=x_{j,n}$ and
$$    \beta_j(t_1,\ldots,t_{2n}) = \pi_{2n}(x_{j,n}) \lambda_{j,n} ,  $$
where $\lambda_{j,n}$ are the Christoffel numbers corresponding to the
orthogonal polynomial $p_n(x)$ for the weight function \thetag{2.3}.
\endroster
\endproclaim

\demo{Proof}

(1) The orthogonal functions  obtained by orthogonalizing
the system
of basis functions \thetag{1.1} on $[-1,1]$ with weight function $w(x)$ are of
the form
$$    r_{n+1}(t_1,\ldots,t_{n+1};x) = \sum_{k=1}^{n+1} \frac{a_k}{1+t_kx} =
\frac{p_{n}(x)}{\pi_{n+1}(x)} \qquad (a_{n+1} = 1),  $$
where
 $p_{n}(x)$ is a polynomial of degree $n$. The orthogonality implies that
$$   \int_{-1}^1 \frac{p_n(x)}{\pi_{n+1}(x)} \sum_{j=1}^{n}
\frac{c_j}{1+t_jx} \ w(x) \, dx = 0 $$
for every choice of parameters $c_1,c_2,\ldots,c_n$. This is the same as
saying that
$$   \int_{-1}^1 p_n(x)q_{n-1}(x) \frac{w(x)}{\pi_n(x)\pi_{n+1}(x)} \, dx
 = 0 $$
for every polynomial $q_{n-1}(x)$ of degree at most $n-1$. This shows that
$p_n(x)$ is the orthogonal polynomial of degree $n$ for the weight function
\thetag{2.1}.
 The Gauss-Jacobi quadrature for the weight function \thetag{2.1} gives
$$   \int_{-1}^1  \frac{q_{2n-1}(x)w(x)}{\pi_n(x)\pi_{n+1}(x)} \, dx
  = \sum_{j=1}^n \lambda_{j,n} q_{2n-1}(x_{j,n}), \tag 2.4 $$
where $x_{j,n}$ are the zeros of $p_n(x)$ and $q_{2n-1}(x)$ is a polynomial
of degree at most $2n-1$.
If $t_1=0$ then $\pi_{n}(x)$ is a polynomial of degree $n-1$.
In this case
choose $q_{2n-1}(x) = \pi_n(x)\pi_{n+1}(x)/(1+t_kx)$ to obtain
$$   \int_{-1}^1 \frac{w(x)}{1+t_kx}  \,dx =
    \sum_{j=1}^n
      \frac{\lambda_{j,n} \pi_n(x_{j,n})\pi_{n+1}(x_{j,n})}{1+t_kx_{j,n}} ,
 \qquad k=1,2,\ldots,n+1,      $$
from which the assertion follows.

(2)
Let $x_i \equiv x_i(t_1,\ldots,t_{2n})$ and $\beta_i \equiv
\beta_i(t_1,\ldots,t_{2n})$ $(i=1,2,\ldots,n)$. In order to obtain
the nodes $x_i$ and weights $\beta_i$, one needs to solve the system of
equations
$$    \int_{-1}^1 \frac{w(x)}{1+t_jx}\, dx = \sum_{i=1}^n
\frac{\beta_i}{1+t_jx_i}\ , \qquad j=1,2,\ldots,2n. $$
This is equivalent to finding $x_i$ and $\beta_i$ $(i=1,2,\ldots,n)$ such that
$$  \int_{-1}^1 \frac{q_{2n-1}(x)}{\pi_{2n}(x)}w(x) \, dx
   = \sum_{i=1}^n \beta_i \frac{q_{2n-1}(x_i)}{\pi_{2n}(x_i)} \tag 2.5 $$
for every polynomial $q_{2n-1}$ of degree at most $2n-1$. From the theory
of orthogonal polynomials (Gauss-Jacobi quadrature) we then know that
the nodes $x_i$ are the zeros of the orthogonal polynomial of degree $n$
for the weight function \thetag{2.3} and then
$$  \int_{-1}^1 \frac{q_{2n-1}(x)}{\pi_{2n}(x)}w(x) \, dx
   = \sum_{i=1}^n \lambda_{i,n} q_{2n-1}(x_i) .  $$
Choose $q_{2n-1}(x) = \pi_{2n}(x)/(1+t_kx)$ $(k=1,\ldots,2n)$ to obtain
the assertion. \qed
\enddemo

In case the weight function $w(x)$ is equal to either
$$  w(x) = \cases
    \frac{1}{\sqrt{1-x^2}}, \\
    \sqrt{1-x^2}, \\
    \sqrt{\frac{1-x}{1+x}}, \endcases  $$
the orthogonal polynomials for the modified weight functions
$w_n(x)$ given by \thetag{2.1} or \thetag{2.3} can be found explicitly and
are known as Bernstein-Szeg\H o polynomials \cite{26, p.\ 31}.

The nodes of both, orthogonal quadrature and Gaussian quadrature, are
equal to the zeros of the orthogonal
polynomial $p_{n}(x)$ with respect to a weight function depending on $n$.
Therefore, the distribution of the quadrature
nodes for orthogonal quadrature and Gaussian quadrature is given by the
distribution of the zeros
of orthogonal polynomials with a weight function depending on the degree $n$.
The asymptotic zero distribution and various asymptotic results for
such orthogonal polynomials have been
investigated recently by Mhaskar and Saff \cite{22}--\cite{24}, Gonchar, L\'opez
and Rakhmanov \cite{11} \cite{21} and L\'opez \cite{16} \cite{17} \cite{20}.
An almost immediate corollary to their results is:

\proclaim{Theorem 2}
Suppose that the asymptotic distribution of the parameters $\{t_i: i \leq n\}$
is given by a measure $\nu$ on $[-1,1]$, i.e., for every continuous function
$f$ on $[-1,1]$ we have
$$   \lim_{n \to \infty} \frac1n \sum_{j=1}^n f(t_j) = \int_{-1}^1 f(x)\,
       d\nu(x).    \tag 2.6 $$
If $\nu = p \delta_{-1} + q \delta_{1} + r \delta_0 + [1-(p+q+r)] \nu_0$,
where $\delta_a$ is the unit measure with all its mass concentrated at $a$,
and $p,q,r>0$, $p+q+r \leq 1$ and
$$  \int_{-1}^1 \log |t| \, d\nu_0(t) < \infty , \tag 2.7 $$
then the asymptotic
distribution of the nodes for both, orthogonal and
Gaussian quadrature, is given by the measure
$$   \mu = p\delta_{1} + q\delta_{-1} + r \mu_0 + [1-(p+q+r)]\mu_b,  \tag 2.8 $$
where $\mu_0$ is the arcsin measure on $[-1,1]$ with weight function
$$  \mu_0'(x) = \frac{1}{\pi} \frac{1}{\sqrt{1-x^2}} $$
and $\mu_b$ is an absolutely continuous measure on $[-1,1]$ with weight function
$$   \mu_b'(x) = \frac1{\pi} \frac{1}{\sqrt{1-x^2}} \int_{-1}^1
\frac{\sqrt{1-t^2}}{1+xt} \, d\nu_0(t) .   \tag 2.9 $$
\endproclaim

\demo{Proof}
For orthogonal quadrature the weight function for the orthogonal polynomials
is $w_n(x) = w(x) e^{-2n\varphi_n(x)}$ with
$$  \varphi_n(x) = \frac{1}{2n} \left[\log \pi_{n+1}(x) +
       \log \pi_{n}(x)\right] , $$
whereas for Gaussian quadrature it is
$w_n(x) = w(x) e^{-2n\psi_n(x)}$ with
$$   \psi_n = \frac{1}{2n} \log \pi_{2n}(x) . $$
By the weak convergence \thetag{2.6} we have
$$ \multline
\lim_{n \to \infty} \varphi_n(x) = \lim_{n \to \infty} \psi_n(x) \\
= p \log(1-x) + q \log(1+x) + [1-(p+q+r)] \int_{-1}^1 \log (1+tx)\,
d\nu_0(t), \endmultline \tag 2.10 $$
uniformly on closed sets of $(-1,1)$.
By \cite{11, Theorem on p.\ 124}, the asymptotic distribution
of the zeros of the corresponding orthogonal polynomials
is given by a measure $\mu$ on $[-1,1]$  which is the unique
solution of the integral equation
$$   \int_{-1}^1 \log \frac1{|x-t|} \, d\mu(t) + \varphi(x) = C,
\qquad x \in [-1,1],    \tag 2.11 $$
where $\varphi(x)$  is the right-hand side of \thetag{2.10} and $C$ is a
constant. We easily find
$$ \multline
\varphi(x) = p\log(1-x) + q\log(1+x) - [1-(p+q+r)] \int
\log\frac{1}{|x-t|} \, d\nu_0(-1/t)\\
 + [1-(p+q+r)] \int_{-1}^1 \log|t| \, d\nu_0(t). \endmultline $$
The last term in this expression is a finite constant by \thetag{2.7} and
may be absorbed by the constant in \thetag{2.11}. Notice that $\varphi(x)$
may become $-\infty$ near $\pm 1$, which
at first sight is not allowed when applying the result of \cite{11}, but
a careful analysis of the proof in \cite{11} shows that the case under
consideration is also covered.
The integral
$$  \varphi_b(x) = \int_{(-\infty,-1)\cup(1,\infty)}
\log\frac{1}{|x-t|} \, d\nu_0(-1/t) $$
is the logarithmic potential of the mass distribution $\nu_0(-1/t)$ outside
$[-1,1]$. The integral equation
$$  \int_{-1}^1 \log\frac{1}{|x-t|}\, d\mu_b(t) - \varphi_b(x) = C, \qquad
     x \in [-1,1], $$
has a unique solution $\mu_b$ which is the `balayage' of the mass distribution
$\nu_0(-1/t^2)$ outside $[-1,1]$ to the interval $[-1,1]$ (see, e.g., \cite{21,
pp.\ 126--127}), and this measure $\mu_b$ is absolutely continuous with
weight function given by \thetag{2.9}. The arcsin distribution has a constant
logarithmic potential on $[-1,1]$, which can be absorbed into the constant
$C$ in \thetag{2.11}. Hence, the unique solution of \thetag{2.11}
is the measure given by \thetag{2.8}. \qed
\enddemo

Let us take a closer look at the limiting zero distribution measure $\mu$ given
by \thetag{2.8}. If the measure $\nu$ has all its mass at the point $a \in
(-1,1)$ then
$$   \mu'(x) = \frac1{\pi} \frac{1}{\sqrt{1-x^2}}
\frac{\sqrt{1-a^2}}{1+ax} .   \tag 2.12 $$
If $a=0$, then the limiting zero distribution is the well-known arcsin
distribution. For $a=-1$ the limiting zero distribution is a degenerate
distribution at $1$ wheras for $a=1$ the limiting zero distribution
is concentrated at $-1$.
If $a$ is close to $1$ in \thetag{2.12}, then $\mu$ has even more mass
concentrated around $-1$ than the arcsin distribution.
The measures with weight functions given by \thetag{2.9} or
\thetag{2.12}
are thus somewhere between an arcsin measure and  degenerate measures at
$\pm 1$.

Even though the asymptotic distribution of the nodes for orthogonal quadrature
and Gaussian quadrature is the same, we see that the $n$ nodes for the Gaussian
quadrature formula depend on $2n$ parameters $t_i$ whereas the $n$ nodes
for orthogonal quadrature depend on $n+1$ parameters $t_i$. The nodes for
Gaussian quadrature therefore use more (asymptotic) information of the basis
functions.

\head 3. Properties of the nodes \endhead

The orthogonal functions $r_n(t_1,\ldots,t_n;x)$ $(n=1,2,\ldots)$ are an
extended Haar system on $[-1,1]$, and therefore the properties of the zeros
given by Pinkus, Ziegler and Videnski\u\i\  are valid. For this special Haar
system additional useful properties can be proved:

\proclaim{Theorem 3a}
\roster
\item The zeros of $r_n(t_1,\ldots,t_n;x)$ and
$r_{n+1}(t_1,\ldots,t_n,t_{n+1};x)$
on $[-1,1]$ are strictly interlacing.
\item The zeros of $r_n(t_1,\ldots,t_{n-1},t_n;x)$ and $r_n(t_1,\ldots,t_{n-1},
t_{n+1};x)$ on $[-1,1]$ are strictly interlacing.
\item The zeros of $r_n(t_1,\ldots,t_{i-1},t_{i+1},\ldots,t_n,t_i;x)$
$(1 \leq i \leq n)$ and $r_n(t_1,\ldots,\allowbreak t_n;x)$ are strictly
interlacing. \endroster
\endproclaim

\demo{Proof} Property (1) is proved in Pinkus and Ziegler \cite{25} and
Videnski\u\i\  \cite{28}, and property (2) can be found in Pinkus and
Ziegler \cite{25}. For property (3) we notice that
$$ \align
r_n(t_1,\ldots,t_n;x) = &r_{n-1}(t_1,\ldots,t_{i-1},t_{i+1},\ldots,t_n;x) \\
     &- \kappa_1 r_{n-1}(t_1,\ldots,t_{i-1},t_{i+1},\ldots,t_{n-1},t_i;x)
\endalign     $$
and
$$ \align
r_n(t_1,\ldots,t_{i-1},t_{i+1},\ldots,t_n,t_i;x) =
&r_{n-1}(t_1,\ldots,t_{i-1},t_{i+1},\ldots,t_{n-1},t_i;x) \\
     &- \kappa_2 r_{n-1}(t_1,\ldots,t_{i-1},t_{i+1},\ldots,t_{n-1},t_n;x),
\endalign     $$
where
$$  \kappa_1 = \frac{\left( r_{n-1}(t_1,\ldots,t_{i-1},t_{i+1},\ldots,t_n;x)
, \frac{1}{1+t_ix} \right)}{\left\|
r_{n-1}(t_1,\ldots,t_{i-1},t_{i+1},\ldots,t_{n-1},t_i;x) \right\|_2^2} $$
and
$$  \kappa_2 = \frac{\left(
r_{n-1}(t_1,\ldots,t_{i-1},t_{i+1},\ldots,t_{n-1},t_i;x)
, \frac{1}{1+t_nx} \right)}{\left\|
r_{n-1}(t_1,\ldots,t_{i-1},t_{i+1},\ldots,t_n;x) \right\|_2^2}. $$
Let $x_i$ $(i=1,2,\ldots,n-1)$ be the zeros of $r_n(t_1,\ldots,t_n;x)$; then
$$ \multline
r_n(t_1,\ldots,t_{i-1},t_{i+1},\ldots,t_n,t_i;x_j)  \\
    = r_{n-1}(t_1,\ldots,t_{i-1},t_{i+1},\ldots,t_{n-1},t_i;x_j) (1-\kappa_1
\kappa_2).
\endmultline \tag 3.1 $$
If we compute the $L_2$-norm of $r_n(t_1,\ldots,t_n;x)$, then
$$ \align
\| r_n(t_1,\ldots,t_n;x) \|_2^2 &= \left( r_n(t_1,\ldots,t_n;x),
\frac{1}{1+t_nx} \right)  \\
   &= \left( r_{n-1}(t_1,\ldots,t_{i-1},t_{i+1},\ldots,t_n;x), \frac{1}{1+t_nx}
\right) \\
& \quad - \kappa_1 \left(
r_{n-1}(t_1,\ldots,t_{i-1},t_{i+1},\ldots,t_{n-1},t_i;x) , \frac{1}{1+t_nx}
\right) \\
  &= \| r_{n-1}(t_1,\ldots,t_{i-1},t_{i+1},\ldots,t_n;x) \|_2^2 (1 - \kappa_1
\kappa_2)
 \endalign $$
so that $1-\kappa_1 \kappa_2 > 0$. Then property (3) follows immediately
from \thetag{3.1} and property (1). \qed
\enddemo

We have similar results for the nodes corresponding to Gaussian quadrature.
Now the order of the parameters $t_1,\ldots,t_{2n}$ is irrelevant.

\proclaim{Theorem 3b}
Let $x_i(t_1,\ldots,t_{2n})$ $(i=1,2,\ldots,n)$ be
the quadrature nodes for maximum accuracy.
\roster
\item The $n$ nodes $x_i(t_1,\ldots,t_{2n})$ and the $n+1$ nodes
$x_i(t_1,\ldots,t_{2n+1},t_{2n+2})$ are strictly interlacing.
\item The $n$ nodes $x_i(t_1,\ldots,t_{2n-1},t_{2n})$ and the $n$ nodes
$x_i(t_1,\ldots,t_{2n-1},t_{2n+1})$ are strictly interlacing.
\endroster
\endproclaim

\demo{Proof of (1)}
Let $x_0(t_1,\ldots,t_{2n}) = -1$ and $x_{n+1}(t_1,\ldots,t_{2n}) = 1$. \newline
Suppose that the nodes $x_i(t_1,\ldots,t_{2n})$ and
$x_i(t_1,\ldots,t_{2n+2})$ are not strictly interlacing; then either
\newline
{\bf (a)}
$x_i(t_1,\ldots,t_{2n}) = x_j(t_1,\ldots,t_{2n+2})$ for some $i$ and $j$
with $1 \leq i \leq n$ and $1 \leq j \leq n+1$ or \newline
{\bf (b)}
$x_j(t_1,\ldots,t_{2n+2}), x_{j+1}(t_1,\ldots,t_{2n+2}) \in \left(
(x_i(t_1,\ldots,t_{2n}),x_{i+1}(t_1,\ldots,t_{2n}) \right)$ for some $i$ and $j$
with $0 \leq i \leq n$ and $1 \leq j \leq n$. \newline
We will show that both (a) and (b) are impossible. \newline
{\bf Case (a):}
Suppose that there are $\ell$ nodes $x_i(t_1,\ldots,t_{2n+2})$, $i \in
\{i_1,i_2,\ldots,i_{\ell}\}$, such that $x_i(t_1,\ldots,t_{2n+2})
\in \{x_j(t_1,\ldots,t_{2n}): 1 \leq j \leq n\}$; then obviously $\ell \leq n$,
and thus there exists a $k$ such that $x_k(t_1,\ldots,t_{2n+2}) \notin
\{x_j(t_1,\ldots,t_{2n}): 1 \leq j \leq n\}$.
Let $\phi$ be a linear combination of the first $2n$ basis functions given in
\thetag{1.1} such that $\phi\left( x_i(t_1,\ldots,t_{2n})\right) = 0$
for $i = 1,2,\ldots,n$ and $\phi\left( x_i(t_1,\ldots,t_{2n+2}) \right)=0$ for
$i=1,2,\ldots,2n+1$ but $i \neq k$. Then we already have $2n-\ell$ equations
to determine $\phi$. In addition we add the equations
$\phi(z_i) = 0$ for $i=1,\ldots,\ell-1$, where the $z_i$ differ from the
$2n-\ell$ points that we already used. For such $\phi$ we have
$$   \int_{-1}^1 \phi(x)w(x)\, dx = \sum_{i=1}^n \beta_i(t_1,\ldots,t_{2n})
   \phi\left(x_i(t_1,\ldots,t_{2n})\right) = 0. $$
But we also have
$$  \align
\int_{-1}^1 \phi(x)w(x)\, dx &= \sum_{i=1}^{n+1} \beta_i(t_1,\ldots,t_{2n+2})
\phi\left(x_i(t_1,\ldots,t_{2n+2})\right) \\
     &= \beta_k(t_1,\ldots,t_{2n+2})
\phi\left(x_k(t_1,\ldots,t_{2n+2})\right) \neq 0 ,
   \endalign $$
because $\phi$ can have at most $2n-1$ zeros and $\beta_k > 0$. This gives a
contradiction. \newline
{\bf Case (b):}
Let $\phi$ be the linear combination of the first $2n$ basis functions given
in  \thetag{1.1} such that $\phi\left( x_i(t_1,\ldots,t_{2n})\right) = 0$ for
$i=1,2,\ldots,n$, and $\phi\left( x_i(t_1,\ldots,t_{2n+2}) \right) = 0$ for
$i=1,2,\ldots,n+1$ but with $i \neq j$ and $i \neq j+1$. This gives $2n-1$
conditions to determine $\phi$, and for such $\phi$ we have
$$   \int_{-1}^1 \phi(x)w(x)\, dx = \sum_{i=1}^n \beta_i(t_1,\ldots,t_{2n})
   \phi\left(x_i(t_1,\ldots,t_{2n})\right) = 0. $$
On the other hand we also have
$$  \align
\int_{-1}^1 \phi(x)w(x)\, dx &= \sum_{i=1}^{n+1} \beta_i(t_1,\ldots,t_{2n+2})
\phi\left(x_i(t_1,\ldots,t_{2n+2})\right) \\
     &= \beta_j(t_1,\ldots,t_{2n+2})
\phi\left(x_j(t_1,\ldots,t_{2n+2})\right) \\
&+ \quad \beta_{j+1}(t_1,\ldots,t_{2n+2})
\phi\left(x_{j+1}(t_1,\ldots,t_{2n+2})\right).
   \endalign $$
By assumption we have sign $\phi\left(x_j(t_1,\ldots,t_{2n+2})\right) =
\text{sign }\phi\left(x_{j+1}(t_1,\ldots,t_{2n+2})\right)$, and since the
weights $\beta_j$ and $\beta_{j+1}$ are positive, we have a contradiction.
\enddemo
\demo{Proof of (2)}
Suppose that $x_i\equiv x_i(t_1,\ldots,t_{2n})$ and
$y_i\equiv x_i(t_1,\ldots,t_{2n-1},t_{2n+1})$ are not strictly interlacing;
then either \newline
{\bf (a)}
$x_i = y_i$ for all $i=1,2,\ldots,n$; or \newline
{\bf (b)}
$x_i = y_j$ for some $i$ and $j$ with $1 \leq i \leq n$ and $1 \leq j
\leq n$; or \newline
{\bf (c)}
$y_j,y_{j+1} \in (x_i,x_{i+1})$ for some $i$ and $j$ with $0 \leq i \leq n$ and
$1 \leq j \leq n-1$.\newline
{\bf Case (a):}
The system of basis functions \thetag{1.1} is a Haar system on $[-1,1]$, hence
we can find a function $f_k$ in the linear space spannen by the first
$2n-1$ basis function in \thetag{1.1} such that
$$ \align
f_k(x_i) &= \delta_{i,k}, \qquad i=1,2,\ldots,n , \\
f_k(\xi_i) &= 0, \qquad i=1,2,\ldots, n-1
\endalign $$
where $\xi_i \in (-1,1) \setminus \{x_1,x_2,\ldots,x_n\}$ are $n-1$
arbitrary points. Then on one hand
$$    \int_{-1}^1 f_k(x) w(x) \, dx = \sum_{i=1}^n \beta_i(t_1,\ldots,t_{2n})
  f_k(x_i) = \beta_k(t_1,\ldots,t_{2n}) , $$
and on the other hand (since $x_i=y_i$ for $1 \leq i \leq n$)
$$    \int_{-1}^1 f_k(x) w(x) \, dx =
\sum_{i=1}^n \beta_i(t_1,\ldots,t_{2n-1},t_{2n+1})
  f_k(x_i) = \beta_k(t_1,\ldots,t_{2n-1},t_{2n+1}) . $$
Hence
$\beta_i(t_1,\ldots,t_{2n-1},t_{2n+1}) = \beta_i(t_1,\ldots,t_{2n})$ for
$i=1,2,\ldots,n$. Therefore the quadrature formula will give the exact
result for every function $\frac{1}{1+t_ix}$ with $i=1,2,\ldots,2n+1$.
Let $\phi$ be a linear combination of these $2n+1$ basis functions such
that $\phi(x_i)=0$ and $\phi'(x_i)=0$ $(i=1,2,\ldots,n)$; then
$\phi$ does not changes sign on $[-1,1]$ because $\phi(x)
= \pm p_n^2(x)/\pi_{2n+1}(x)$ with $p_n$ a polynomial of degree $n$ and
$\pi_{2n+1}$ given by \thetag{2.2}. This implies that
$$   \int_{-1}^1 \phi(x)w(x) \, dx \neq 0 . $$
This contradicts
$$     \int_{-1}^1 \phi(x)w(x) \, dx = \sum_{i=1}^n \beta_i \phi(x_i) = 0 . $$
{\bf Cases (b) and (c):}
one can obtain a contradiction in a way similar to cases (a) and (b) of (1).
\qed
\enddemo

The following result of A.A. Markov (see, e.g., \cite{26, Thm.\ 6.12.1 on p.\
115}) allows us to give more information on the behavior of the nodes
$x_i(t_1,\ldots,t_n,t)$ for orthogonal quadrature, and
$x_i(t_1,\ldots,t_{2n-1},t)$ for Gaussian quadrature, as functions of $t \in
(-1,1)$:

\proclaim{Theorem} {\rm (Markov)}
Let $w(x,t)$ be a weight function on $[a,b]$ depending on a parameter $t$
such that $w(x,t)$ is positive and continuous for $a < x < b$ and $c < t <
d$. Suppose that $w_t(x,t) = \frac{\partial}{\partial t} w(x,t)$ exists and is
continuous for $a < x < b$ and $c < t < d$, and that the moments
$$   \int_a^b x^k w_t(x,t) \, dx , \qquad  k=0,1,2,\ldots,2n-1 , $$
converge uniformly in every closed subinterval $c' \leq t \leq d'$ of
$(c,d)$. If the zeros of the orthogonal polynomial $p_n(x,t)$ with weight
function $w(x,t)$ are denoted by $x_i(t)$ $(i=1,2,\ldots,n)$ and if
$w_t(x,t)/w(x,t)$ is a strictly increasing (decreasing) function of $x \in
(a,b)$, then
$x_i(t)$ is a continuously differentiable and
strictly increasing (decreasing) function of $t \in (c,d)$ for
every fixed $i$. \endproclaim

\proclaim{Theorem 4} \roster
\item If $x_i(t_1,\ldots,t_n,t)$ is the $i$th node for  orthogonal
quadrature with $n$ nodes,
then $x_i(t_1,\ldots,t_n,t)$ is a continuously differentiable and
strictly decreasing function of $t \in (-1,1)$.
\item
Every node $x_i(t_1,\ldots,t_{2n-1},t)$ for Gaussian quadrature with $n$ nodes
is a continuously differentiable and strictly decreasing function of $t \in
(-1,1)$. \endroster
\endproclaim

\demo{Proof}
By Theorem 1 (1)
we know that $x_i(t_1,\ldots,t_n,t)$ is the $i$th zero
of the orthogonal polynomial $p_{n}(x,t)$ of degree $n$ for the weight
function
$$    w(x,t) = \frac{w(x)}{\pi_{n}^2(x)(1+tx)} . $$
A simple computation gives
$$   \frac{w_t(x,t)}{w(x,t)} = - \frac{x}{1+tx} $$
and this is a decreasing function of $x \in (-1,1)$. The result thus follows
by Markov's theorem.

Similarly, by Theorem 1 (2), for the nodes
$x_i(t_1,\ldots,t_{2n-1},t)$ of Gaussian quadrature
one can use Mar\-kov's theorem  for the weight function
$$   w(x,t) = \frac{w(x)}{\pi_{2n-1}(x)(1+tx)},$$
to obtain the desired result. \qed
\enddemo

We can also give some monotonicity results for the first and last
quadrature coefficients. Let
$\beta_i(t_1,\ldots,t_n,t)$ be the quadrature coefficients
corresponding to the nodes $x_i(t_1,\ldots,t_n,t)$ for
orthogonal quadrature with $n$ nodes, and let $t_1 = 0$; then
the first quadrature coefficient
$\beta_1(t_1,\ldots,t_n,t)$ is a strictly decreasing continuously differentiable
function of $t \in (-1,1)$, and the last quadrature coefficient
$\beta_n(t_1,\ldots,t_n,t)$ is a strictly
increasing continuously differentiable function of $t \in (-1,1)$.
Similarly, if $\beta_i(t_1,\ldots,t_{2n-1},t)$ are the quadrature coefficients
for Gaussian
quadrature with $n$ nodes (with $t_1=0$); then $\beta_1(t_1,\ldots,t_{2n-1},t)$
is a strictly decreasing continuously differentiable function of $t \in (-1,1)$
and $\beta_n(t_1,\ldots,t_{2n-1},t)$ is a strictly increasing continuously
differentiable function of $t \in (-1,1)$. The proofs of these results
are rather technical (see \cite{27}) and since only the
extreme quadrature weights are considered we decided not to include them here.

\head 4. Convergence results \endhead

If the integral
$$ \int_{-1}^1 f(x)w(x)\,dx $$
is approximated by the orthogonal quadrature sum with $n$ nodes,
$$  \sum_{j=1}^{n} \beta_i(t_1,\ldots,t_{n+1}) f(x_i(t_1,\ldots,t_{n+1})), $$
we denote the error by
$$   E_n^o f = \int_{-1}^1 f(x)w(x)\,dx  - \sum_{j=1}^{n}
\beta_i(t_1,\ldots,t_{n+1}) f(x_i(t_1,\ldots,t_{n+1})).  $$
The error of the Gaussian quadrature formula is denoted by
$$   E_n^g f = \int_{-1}^1 f(x)w(x)\,dx  - \sum_{j=1}^n
\beta_i(t_1,\ldots,t_{2n}) f(x_i(t_1,\ldots,t_{2n})).  $$
It is clear that these quadrature formulas are only going to be relevant
if for a large enough class of functions $f$ one has
$E_n^o f \rightarrow 0$ as $n \rightarrow \infty$ and/or
$E_n^g f \rightarrow 0$ as $n \rightarrow \infty$.
It is well known that the sequence of basis functions
\thetag{1.1} is dense in $C[-1,1]$ with respect to the uniform norm if and
only if
$$        \sum_{k=1}^\infty (1 - |c_k|) = \infty,    \tag 4.1 $$
where $c_k = -1/t_k - \sqrt{1/t_k^2-1}$ \cite{1, p.\ 254, \S 7}\footnote{There
is a misprint in this reference: `interior' should be replaced by `exterior'.
We thank P. Borwein for pointing out this misprint}. The following result
is straightforward, but we mention it for the sake of completeness.

\proclaim{Theorem 5}
Let $t_1,t_2,t_3,\ldots$ be distinct parameters in $(-1,1)$ and
 $\phi_{n+1}(x)$ be the best approximation of $f \in C[-1,1]$ using
linear combinations of
$\frac{1}{1+t_1x},\frac{1}{1+t_2x},\ldots,\allowbreak \frac{1}{1+t_{n+1}x}$.
\roster
\item If $t_1=0$, then for orthogonal quadrature, there exists a positive
constant $M$ not depending on $n$ such that
$$    |E_n^o f| \leq M \| f- \phi_{n+1} \|_{\infty}. $$
As a consequence, $E_n^o f \rightarrow 0$ as $n \rightarrow \infty$ for every
bounded Riemann integrable function on $[-1,1]$ whenever \thetag{4.1} holds.
\item  For Gaussian quadrature we have similarly
$$    |E_n^g f| \leq M \| f- \phi_{2n} \|_{\infty}. $$
As a consequence, $E_n^g f \rightarrow 0$ as $n \rightarrow \infty$ for every
bounded Riemann integrable function on $[-1,1]$ whenever \thetag{4.1} holds.
\endroster
\endproclaim

\demo{Proof}
The result for orthogonal quadrature
 can be proved as in the case of polynomial
interpolation \cite{8, pp.\ 126--129}.
For Gaussian quadrature we only need to show that
$\sum_{i=1}^n \beta_i(t_1,\ldots,t_{2n})$ is bounded. If $t_1=0$, then
$$  \sum_{i=1}^n \beta_i(t_1,\ldots,t_{2n})  = \int_{-1}^1 w(x)\, dx. $$
If $t_1 \neq 0$, then
$$  \int_{-1}^1 \frac{w(x)}{1+t_1x} \, dx = \sum_{i=1}^n
\frac{\beta_i}{1+t_1x_i} . $$
The function $1/(1+t_1x)$ has the minimum value $1/(1+|t_1|)$ on $[-1,1]$, hence
$$   \frac{1}{1+|t_1|} \sum_{i=1}^n \beta_i(t_1,\ldots,t_{2n})
   < \int_{-1}^1 \frac{w(x)}{1+t_1x}\, dx, $$
which shows the desired boundedness. \qed
\enddemo

Let $f \in C[-1,1]$ and denote by $L_nf$ the rational interpolant of $f$ at the
points $x_i(t_1,\ldots,t_{n+1})$ $(i=1,2,\ldots,n)$ using the basis functions
$\frac{1}{1+t_1x},\frac{1}{1+t_2x},\ldots,\frac{1}{1+t_{n}x}$. The previous
theorem implies that
$$   \lim_{n \rightarrow \infty} \int_{-1}^1 (L_nf)(x) w(x) \, dx
   = \int_{-1}^1 f(x) w(x) \, dx , $$
whenever \thetag{4.1} holds.
The rate of convergence depends on how well $f$ can be approximated in
the uniform norm by the basis functions \thetag{1.1}, and this depends on the
choice of the parameters $t_i$. We now show that rational interpolation
converges in $L_2(w)$ and thus also in $L_1(w)$, thereby extending the famous
Erd\H os-Tur\'an result for polynomial interpolation \cite{6}.

\proclaim{Theorem 6}
Suppose $t_1=0$ and $t_2,t_3,\ldots$ are distinct parameters in $(-1,1)$ such
that \thetag{4.1} holds.
If $f \in C[-1,1]$ and if $L_nf(x) = \sum_{i=1}^{n} \frac{a_i}{1+t_ix}$
interpolates
$f$ at the zeros $x_i(t_1,\ldots,t_{n+1})$ of $r_{n+1}(t_1,\ldots,t_{n+1};x)$,
then $\|L_nf - f\|_2 \rightarrow 0$ as $n \rightarrow \infty$.
\endproclaim

\demo{Proof}
Let $\phi_{n}(x)$ be the best uniform approximation of $f$ using the
basis functions
$\{\frac{1}{1+t_1x},\ldots,\frac{1}{1+t_{n}x}\}$; then
$$  \| L_nf - f \|_2 \leq \| f - \phi_{n} \|_2 + \| L_nf - \phi_{n} \|_2 .
$$
Clearly, $\| f - \phi_{n} \|_2 \rightarrow 0$ as $n \rightarrow \infty$.
Now $[(L_nf)(x) - \phi_{n}(x)]^2$ is a linear combination of the functions
$x^j/[\pi_{n}^2(x)]$ $(j=0,1,\ldots,2n-2)$, and by \thetag{2.4} the
quadrature formula is therefore exact for $[(L_nf)(x) - \phi_{n}(x)]^2$. This
gives
$$ \align
\|L_nf - \phi_{n}\|_2^2 &= \sum_{i=1}^{n} \beta_i [(L_nf)(x_i) -
    \phi_{n}(x_i)]^2  \\
   &= \sum_{i=1}^{n} \beta_i [f(x_i) - \phi_{n}(x_i)]^2 \\
   & \leq \|f-\phi_{n} \|_\infty^2 \int_{-1}^1 w(x)\, dx .
 \endalign $$
The result now follows because $\|f-\phi_{n}\|_\infty \rightarrow 0$ as
$n \rightarrow \infty$. \qed
\enddemo

\head 5. Generalizations \endhead

All the properties and results for orthogonal quadrature that we have given
above are also true for
the orthogonal system obtained by orthogonalizing the system of basis functions
$$ \frac{1}{1+t_ix}, \left. \frac{\partial}{\partial t} \frac{1}{1+tx}
   \right|_{t=t_i}, \ldots \left. \frac{\partial^{\mu_i-1}}{\partial
t^{\mu_i-1}} \frac{1}{1+tx} \right|_{t=t_i}, \qquad i=1,\ldots,k,  \tag 5.1 $$
where $k \in \Bbb N$, $\mu_i \in \Bbb N$ and $\sum_{i=1}^k \mu_i = n+1$.
If the parameters $t_i$ are all distinct and in $(-1,1)$, then this system of
functions is an extended Haar system on $[-1,1]$. If $t \notin
\{t_1,\ldots,t_k\}$ then
$$    r_{n+1}(t_1^{(\mu_1)},\ldots,t_k^{(\mu_k)},t;x) = \frac{1}{1+tx}
 + \sum_{i=1}^k  \sum_{j=0}^{\mu_i-1} a_{i,j} \left.
\frac{\partial^j}{\partial t^j} \frac{1}{1+tx} \right|_{t=t_i}  $$
is orthogonal to the system \thetag{5.1} with respect to the weight function
$w(x)$; if $t=t_\ell$ for some $\ell$ with $1 \leq \ell \leq k$ then
$$ \gather
r_{n+1}(t_1^{(\mu_1)},\ldots,t_{\ell-1}^{(\mu_{\ell-1})},
t_{\ell+1}^{(\mu_{\ell+1})},\ldots,t_k^{(\mu_k)},t_{\ell}^{(\mu_{\ell}+1)};x)
 \qquad \qquad  \\
 = \left. \frac{\partial^{\mu_\ell}}{\partial t^{\mu_\ell}} \frac{1}{1+tx}
\right|_{t=t_\ell} + \sum_{i=1}^k  \sum_{j=0}^{\mu_i-1} a_{i,j} \left.
\frac{\partial^j}{\partial t^j} \frac{1}{1+tx} \right|_{t=t_i}
  \endgather $$
is orthogonal to the system \thetag{5.1} with respect to the weight function
$w(x)$. Because of the Haar property it follows that these orthogonal functions
$r_{n+1}$ have precisely $n$ zeros in $[-1,1]$ which are all simple
\cite{25}, \cite{28}. These orthogonal functions share the same properties as
the orthogonal functions for the system \thetag{1.1}:  the reason for this
is that Theorem 1 can be extended to this system.
Similarly all the results and properties for Gaussian
quadrature that we have given are also true for the system of basis functions
\thetag{5.1} with $k \in {\Bbb N}, \mu_i \in \Bbb N$ and
$\sum_{i=1}^k \mu_i = 2n$.
The only things that change are the weight functions \thetag{2.1} and
\thetag{2.3},
which need to be appropriately modified.

\proclaim{Theorem 7} \roster
\item The zeros of $r_{n+1}(t_1^{(\mu_1)},\ldots,t_k^{(\mu_k)};x)$ are also the
zeros of the orthogonal polynomial $p_{n}(x)$ of degree $n$ for
the weight function
$$   \frac{w(x)}{\prod_{i=1}^{k-1} (1+t_ix)^{2\mu_i} (1+t_kx)^{2\mu_k-1}}. $$
\item The $n$ nodes for Gaussian quadrature with the system \thetag{5.1}
are the zeros of the orthogonal polynomial $p_n(x)$ of degree $n$ for the weight
function
$$     \frac{w(x)}{\prod_{i=1}^k (1+t_ix)^{\mu_i}} . $$
\endroster
\endproclaim
This theorem
allows us to use Markov's theorem to deduce properties of the nodes and the
weights for orthogonal quadrature, but since everything is along the same
lines as before, we do not give the details.

\head 6. Numerical examples \endhead

From Theorem 5 it follows that rational quadrature rules will be effective
whenever the function to be integrated can be approximated well using the
basis functions \thetag{1.1}. Gautschi \cite{9} has shown that functions which
have an infinite number of poles but are regular otherwise, can be integrated
most effectively by Gauss rational quadrature rules by choosing the parameters
$t_i$ in an optimal way. Gautschi \cite{10} explicitly used such quadrature
rules for the computation
of generalized Fermi-Dirac and Bose-Einstein integrals.

For functions with most mass concentrated near the endpoints of
the interval $[-1,1]$, one could choose the parameters $t_i$ in such a
way that more quadrature nodes are situated near $\pm 1$. If we take
the parameters
$$  t_i = 1 - \frac{1}{\sqrt{i}} \qquad i \in {\Bbb N}, \tag 6.1 $$
then there will be more quadrature nodes near $-1$ than near $+1$.
In fact, from Theorem 2 it follows that the asymptotic distribution
of the nodes for both, orthogonal quadrature and Gaussian quadrature,
is degenerate at the point $-1$. Let GR($n$) be the $n$-point Gaussian
quadrature rule, and OR($n$) the $n$-point orthogonal quadrature rule.
The nodes and weights for $n=6$ are given in Table 1. Our computations
have been made on an IBM 3090/400E VF using Brent's multiple precision
package \cite{5}, and an independent computation was made using
{\smc mathematica}$^{TM}$ \cite{29}\footnote{{\smc Mathematica}
is a trademark of Wolfram Research Inc.}.
\vskip 1.5em

\hfil\vbox{\offinterlineskip
\hrule
\halign{&\vrule#&\strut\quad\hfil#\quad\cr
height2pt&\omit&&\omit&\cr
&\hfil nodes for GR(6) \hfil&& weights for GR(6) \hfil&\cr
height2pt&\omit&&\omit&\cr
\noalign{\hrule}
height2pt&\omit&&\omit&\cr
& -0.9797390942708352 &&   0.0528758827013522 &\cr
& -0.8853794251591486 &&   0.1411615118844550 &\cr
& -0.6822351336410264 &&   0.2748067575758422 &\cr
& -0.3156675377072605 &&   0.4657849717765712 &\cr
&  0.2408527285476740 &&   0.6221630733806293 &\cr
&  0.8155273184304977 &&   0.4432078026811501 &\cr
height2pt&\omit&&\omit&\cr
\noalign{\hrule}\cr
height2pt width0pt\cr
\noalign{\hrule}
height2pt&\omit&&\omit&\cr
&\hfil nodes for OR(6) \hfil&& weights for OR(6) \hfil&\cr
height2pt&\omit&&\omit&\cr
\noalign{\hrule}
height2pt&\omit&&\omit&\cr
& -0.9736320979338328 &&   0.0685126325838336 &\cr
& -0.8537169072027923 &&   0.1760476819554412 &\cr
& -0.6094091127142633 &&   0.3192517203251832 &\cr
& -0.2057016948376719 &&   0.4878639628808742 &\cr
&  0.3414560761423378 &&   0.5765658940369015 &\cr
&  0.8474273771128526 &&   0.3717581082177663 &\cr
height2pt&\omit&&\omit&\cr}
\hrule}\hfil
\vskip 0.5em
\centerline{Table 1: nodes and weights for the quadrature rules GR(6) and
OR(6)}
\vskip 0.5em

The nodes and weights for GR($n$) and OR($n$) have been computed by
using the zeros and Christoffel numbers for the orthogonal polynomial
of degree $n$ with the corresponding weight \thetag{2.1} for
orthogonal quadrature and \thetag{2.3} for Gaussian quadrature. These
orthogonal polynomials have been computed using their three-term
recurrence relation. The recurrence coefficients were obtained using
Chebyshev's algorithm with modified moments and the latter were computed
recursively by decomposing the rational weight function into simple
fractions.

Consider the integral
$$   I_1 = \int_{-1}^1 \omega e^{-\omega(x+1)} \, dx, \qquad \omega > 0. $$
This function has more mass near $-1$ than near $+1$, especially for large
$\omega$. In Table~2 we have given the absolute relative errors for the Gaussian
quadrature sum and the orthogonal quadrature sum corresponding to the
rational functions with parameters given by \thetag{6.1},
together with the
Gauss-Legendre quadrature sum GL($n$) with $n$ nodes. The rational quadrature
rules behave in a similar way and give better results than the Gauss-Legendre
quadrature rule.
\vskip 1.5em

\hfil\vbox{\tabskip=0pt \offinterlineskip
\def\tablerule{\noalign{\hrule}}
\halign{\strut#&\vrule#\tabskip=1em plus2em&
  #\hfil&#\hfil&\vrule#&#\hfil&#\hfil&
  #\hfil& \vrule#\tabskip=0pt\cr\tablerule
&& $\omega$ & $n$ && GR($n$) & OR($n$) & GL($n$) &\cr\tablerule
&& 5     & 2  && 0.664(-1) & 0.121     & 0.394     &\cr
&& \omit & 6  && 0.261(-5) & 0.207(-6) & 0.161(-4) &\cr
&& \omit & 10 && 0.215(-8) & 0.842(-10) & 0.212(-5) &\cr
&& \omit & 14 && 0.494(-11) & 0.514(-13) & 0.112(-5) &\cr\tablerule
&& 25    & 2  && 0.978      & 0.988      & 0.999 &\cr
&& \omit & 6  && 0.118(-4)  & 0.255(-4)  & 0.206 &\cr
&& \omit & 10 && 0.397(-8)  & 0.467(-8) & 0.201(-2) &\cr
&& \omit & 14 && 0.206(-12) & 0.122(-12) & 0.153(-4) &\cr\tablerule
}}
\vskip 0.5em
\centerline{Table 2: relative errors for $I_1$}
\vskip 1.5em

The rational quadrature rules perform well because the system of rational
functions with parameters \thetag{6.1} emphasizes the behaviour of the integrand
near $-1$, and even though the quadrature nodes have a degenerate asymptotic
distribution at $-1$, the parameters still satisfy \thetag{4.1} so that
the basis functions are dense in $C[-1,1]$. The fact that the poles
of the rational functions are all close to $-1$ also explains why both
rational quadrature rules give similar results: adding extra rational functions
is not giving a lot of extra information about the integrand.

As a second example we consider the integral
$$    I_2 = \int_{-1}^1 \frac{dx}{\sqrt{(x+3)(x+2)}} . $$
The integrand has a branch cut on the interval $[-3,-2]$ and is in fact
a Stieltjes function given by
$$   \frac{1}{\sqrt{(x+3)(x+2)}} = \frac{1}{\pi} \int_{-3}^{-2} \frac{1}{x-t}
    \frac{dt}{\sqrt{(3+t)(-2-t)}} . $$
Stieltjes functions can be well approximated by rational functions
with poles on the cut. Therefore we will choose the parameters $t_i$
in such a way that the poles of the basis functions \thetag{1.1} are
distributed nicely on $[-3,-2]$. In order to do this we observe that
the zeros of a Chebyshev polynomial of the first kind $T_{3^k}(x)$
are also zeros of Chebyshev polynomials $T_{3^m}(x)$ for every $m \geq k$.
Hence we take $t_1=0$ (so that Theorem 1 can be applied) and
$t_2,t_3,\ldots,t_{3^m+1}$ equal to
$$   t_{i+1} = \frac{-2}{x_{i,3^m} - 5}  , $$
where $x_{i,3^m}$ are the zeros of the Chebyshev polynomial $T_{3^m}(x)$
ordered in such a way that $x_{1,3^m}, x_{2,3^m}, x_{3,3^m}$ are the
zeros of $T_3(x)$ in decreasing order; $x_{4,3^m},\ldots,x_{9,3^m}$ are
those zeros of $T_9(x)$ which are not zeros of $T_3(x)$, again in
decreasing order; and recursively we let $x_{3^{k-1}+1,3^m},\ldots,
x_{3^k,3^m}$ be the zeros of $T_{3^k}(x)$ which are not zeros of
$T_{3^{k-1}}(x)$, in decreasing order. In this way the poles of the
basis functions \thetag{1.1} are dense on $[-3,-2]$ and they are
distributed on $[-3,-2]$ according to the arcsin measure on this interval.
In Table 3 we have given the relative errors for the two rational quadrature
rules and for the Gauss-Legendre rule.
\vskip 1.5em

\hfil\vbox{\tabskip=0pt \offinterlineskip
\def\tablerule{\noalign{\hrule}}
\halign{\strut#&\vrule#\tabskip=1em plus2em&
  #\hfil&\vrule#&#\hfil&#\hfil&
  #\hfil& \vrule#\tabskip=0pt\cr\tablerule
&& $n$ && GR($n$) & OR($n$) & GL($n$) &\cr\tablerule
&&  2  && 0.587(-6) & 0.168(-4) & 0.366(-2) &\cr
&&  6  && 0.161(-19) & 0.262(-16) & 0.567(-7) &\cr
&& 10 && 0.986(-31) & 0.132(-29) & 0.116(-11) &\cr \tablerule
}}
\vskip 0.5em
\centerline{Table 3: relative errors for $I_2$}
\vskip 1.5em

It is clear that the rational quadrature rules give a much better approximation
to the integral $I_2$ than the Gauss-Legendre rule. The reason is obvious: the
integrand is a Stieltjes function which can be approximated much better by
rational functions than by polynomials, especially when we choose the
poles of the rational functions on the interval $[-3,-2]$ on which the Stieltjes
function has its singularities. Also the Gaussian quadrature rule is
slightly better than the orthogonal rule, because the Gaussian rule with $n$
nodes already uses $2n-1$ poles on the interval
$[-3,-2]$, whereas the orthogonal rule only uses $n$ poles.

Finally we consider the integrals
$$   I_3 = \int_{-1}^1  \frac{\pi t/\omega}{\sin (\pi t/\omega)}  \,
     dt, \qquad
   I_4 = \int_{-1}^1 \left( \frac{\pi t/\omega}{\sin (\pi t/\omega)}
    \right)^2 \, dt. $$
Walter Gautschi provided us with numerical results for the approximation
of these two integrals with rational quadrature rules. The integrands contain
infinitely many poles at the points $p_1 = \omega, p_2=-\omega, p_3=2\omega,
p_4= -2\omega, \ldots$. For $I_3$ these are simple poles but for $I_4$
these are poles with multiplicity two. For the orthogonal quadrature rule
we choose the parameters $t_i$ in such a way that they match up with
the poles of the integrand, i.e.,
$$   t_1 = 0, \quad t_{i+1} = -1/p_i, \qquad i=1,2,\ldots. $$
For the Gaussian rule we choose $t_i=-1/p_i$. The results for Gaussian rational
quadrature GR($n$), orthogonal rational quadrature OR($n$) and Gauss-Legendre
quadrature GL($n$) are given in Table 4 and Table 5\footnote{These computations
were done by W. Gautschi on the Cyber 205, with a requested error
tolerance of $0.5\ 10^{-25}$.}.
Both rational quadrature rules perform very well for $I_3$, which is simple to
understand
since these rules explicitly take into account the poles of the integrand.
The Gauss-Legendre rule does not take into account these singularities
of the integrand and as a result it behaves poorly.
\vskip 1.5em

\hfil\vbox{\tabskip=0pt \offinterlineskip
\def\tablerule{\noalign{\hrule}}
\halign{\strut#&\vrule#\tabskip=1em plus2em&
  #\hfil&#\hfil&\vrule#&#\hfil&#\hfil&
  #\hfil& \vrule#\tabskip=0pt\cr\tablerule
&& $\omega$ & $n$ && GR($n$) & OR($n$) & GL($n$) &\cr\tablerule
&& 2     & 5  && 0.46(-11) & 0.26(-8)    & 0.52(-5)  &\cr
&& \omit & 10 && 0.16(-24) & 0.79(-22) & 0.10(-10) &\cr
&& \omit & 15 &&    ---    & 0.17(-24) & 0.20(-16) &\cr\tablerule
&& 1.1   & 5  && 0.21(-8)   & 0.26(-5)   & 0.21(-1) &\cr
&& \omit & 10 && 0.19(-22)  & 0.30(-16)  & 0.26(-3) &\cr
&& \omit & 15 && 0.81(-25)  & 0.67(-25) & 0.32(-5) &\cr\tablerule
&& 1.01   & 5  && 0.43(-8)  & 0.93(-5)   & 0.24 &\cr
&& \omit & 10 && 0.91(-22)  & 0.24(-15)  & 0.64(-1) &\cr
&& \omit & 15 && 0.40(-26)  & 0.61(-26) & 0.17(-1) &\cr\tablerule
}}
\vskip 0.5em
\centerline{Table 4: relative errors for $I_3$}
\vskip 1.5em

When the integrand has poles of second order, then
there is also a difference in the performance of the rational quadrature
rules.

\vskip 1.5em

\hfil\vbox{\tabskip=0pt \offinterlineskip
\def\tablerule{\noalign{\hrule}}
\halign{\strut#&\vrule#\tabskip=1em plus2em&
  #\hfil&#\hfil&\vrule#&#\hfil&#\hfil&
  #\hfil& \vrule#\tabskip=0pt\cr\tablerule
&& $\omega$ & $n$ && GR($n$) & OR($n$) & GL($n$) &\cr\tablerule
&& 2     & 5  && 0.65(-5) & 0.79(-8)    & 0.47(-4)  &\cr
&& \omit & 10 && 0.12(-10) & 0.24(-21) & 0.19(-9) &\cr
&& \omit & 15 && 0.22(-16) & 0.15(-24) & 0.56(-15) &\cr\tablerule
&& 1.1   & 5  && 0.17(-1)   & 0.27(-5)   & 0.13     &\cr
&& \omit & 10 && 0.32(-3)  & 0.30(-16)  & 0.34(-2) &\cr
&& \omit & 15 && 0.38(-5)  & 0.27(-24) & 0.62(-4) &\cr\tablerule
&& 1.01   & 5  && 0.37  & 0.16(-5)   & 0.78 &\cr
&& \omit & 10 && 0.93(-1)  & 0.40(-16)  & 0.40 &\cr
&& \omit & 15 && 0.23(-1)  & 0.11(-24) & 0.16 &\cr\tablerule
}}
\vskip 0.5em
\centerline{Table 5: relative errors for $I_4$}
\vskip 1.5em

For the integral $I_4$
the orthogonal quadrature rule OR($n$) works very well, even for $\omega$
close to one. Both the Gauss-Legendre rule GL($n$) and the rational Gaussian
rule GR($n$) have difficulties when
$\omega$ approaches 1, and even for $\omega = 2$ they do not converge as fast as
the orthogonal quadrature rule OR($n$). For GL($n$) the reason is clear:
no poles of the integrand are being taken into account. On the
other hand GR($n$) approximates the integrand by a rational function
with simple poles, whereas the integrand actually has double poles. Under
estimating the order of the poles thus seems to degrade the convergence.
By choosing $q_{2n-1}(x)=\pi_n(x)\pi_{n+1}(x)/(1+t_kx)^2$ in \thetag{2.4}
we see that the orthogonal quadrature rule integrates
$(1+t_kx)^{-2}$ $(k=2,3,\ldots,n)$ exactly. Hence the orthogonal rule is
more natural when the integrand has poles of order two, but in this
case one could also use the Gaussian quadrature rule corresponding
to the basis function
$$ \frac1{1+t_1x}, \frac1{(1+t_1x)^2} ,
\frac1{1+t_2x}, \frac1{(1+t_2x)^2}
\ldots, \frac1{1+t_nx}, \frac1{(1+t_nx)^2}, $$
which by Theorem 7 uses the zeros of the orthogonal polynomials of degree
$n$ with weight function $\prod_{i=1}^n (1+t_ix)^{-2}$.

\head Acknowledgment \endhead
We are very grateful to Walter Gautschi for his useful comments on an
earlier version of the manuscript and for encouraging us to work out
some numerical examples. We are especially grateful for the numerical
experiments he made in Tables 4 and 5 and for his comments on these
results.

\enddocument